\documentclass{amsart}
\usepackage{amsmath,amssymb}

\newcommand{\R}{\mathbb{R}}
\newcommand{\C}{\mathbb{C}}
\newcommand{\Z}{\mathbb{Z}}

\newcommand{\g}{{\mathfrak{g}}}    
\newcommand{\h}{{\mathfrak{h}}}    

\newcommand{\ord}{\textnormal{ord}}
\newcommand{\Ad}{\textnormal{Ad}}

\let\Box\undefined\newsymbol\Box 1203

\newtheorem{Theorem}{Theorem}[section]
\newtheorem{Lemma}[Theorem]{Lemma}

\newtheorem{Prop}[Theorem]{Proposition}
\newtheorem{Definition}[Theorem]{Definition}
{\theoremstyle{remark}\newtheorem{Remark}[Theorem]{Remark}}
{\theoremstyle{remark}\newtheorem{Example}[Theorem]{Example}}

\begin{document}


%
%
%
\title[Conjugacy classes]{Twisted conjugacy classes, coadjoint
  orbits of loop groups, and D-branes in the WZW model}

\dedicatory{\mbox{\normalfont\mdseries\small \today}}

\author{Stephan Mohrdieck and Robert Wendt} 

\address{Universit\"at Basel\\Rheinsprung 21\\CH--4051 Basle\\Switzerland}
\email{mohrdis@math.unibas.ch}

\address{University of Toronto\\Department of Mathematics\\
    100 St.George Street\\Toronto\\Ontario M5S 3G3\\Canada}
\email{rwendt@math.toronto.edu}

\thanks{We would like to thank A.~Alekseev for helpful discussion and
for pointing out the correct form of the integrality condition to us,
and J.~Fuchs and C.~Schweigert for helpful remarks
  on an earlier version of this paper.
The work of S.M. was supported by the Swiss 
National Science Foundation.}

\begin{abstract}
 We show that untwisted respectively 
twisted conjugacy classes of a compact and simply connected Lie group
which satisfy a certain integrality condition 
correspond naturally to irreducible highest weight representations of
the corresponding affine Lie algebra.
Along the way, we review the classification of  twisted conjugacy classes of a
simply connected compact Lie group $G$ and give a description of their
stabilizers in terms of the Dynkin diagram of the corresponding twisted
affine Lie algebra.
\end{abstract}


\setcounter{page}{1}

\maketitle

\section{Introduction}
Let $G$ be a group and $\tau\in Aut(G)$ an automorphism of $G$. 
The $\tau$-twisted conjugacy classes of $G$ are the orbits 
of the action of $G$ on itself which is given by
$g:h\mapsto gh\tau(g^{-1})$. We shall always allow 
$\tau$ to be the identity in which case the $\tau$-twisted conjugacy 
classes of $G$ are just the ordinary conjugacy classes of $G$.

In this note, we study twisted and untwisted conjugacy classes of a 
compact 
simply connected simple Lie group $G$  which satisfy a certain integrality
condition. This condition appears in the physics
literature where it is used in the classification of 
so called D-branes in the
Wess-Zumino-Witten model  (\cite{AS}, \cite{G}). In a mathematical
context, integral conjugacy classes play a role
in the study of gerbes on compact Lie groups (\cite{Me}). 
Furthermore, the integrality condition for conjugacy 
classes in the group $G$ 
(which we will state in equation (\ref{cond0})
below) can be seen as an
analogue of the integrality condition for coadjoint orbits of 
$G$. It is well known that, via the Borel-Weil-Bott
construction, integral coadjoint orbits of $G$ correspond to
irreducible representations of $G$. On the other hand, it has been
observed in \cite{Me} that untwisted integral conjugacy classes in $G$ are
parametrized by the same set as irreducible highest weight
representations of the affine Lie algebra corresponding to $G$.
This observation has been extended in \cite{St}
to the case of twisted conjugacy classes of $SU(3)$ (where the
twisting automorphism comes from  the diagram automorphism of the Dynkin
diagram of $SU(3)$), and irreducible highest weight representations of
the corresponding twisted affine Lie algebra.

The main goal of this note is to extend the correspondence between
integral conjugacy classes and irreducible highest weight
representations 
to twisted conjugacy classes of arbitrary simple 
compact and simply connected Lie groups and to give a geometric reason
for this correspondence which works for both twisted and untwisted
conjugacy classes simultaneously. 
This is achieved by translating the integrality
condition for untwisted respectively twisted 
conjugacy classes in a compact Lie group $G$
to an integrality condition for coadjoint orbits of the corresponding
untwisted respectively 
twisted loop group of $G$. The correspondence to irreducible highest
weight representation then comes naturally from the ideas of geometric
quantization (or, equivalently, Kirillov's orbit method) which
relates certain coadjoint orbits of a Lie group to its unitary
representations. Along the way, we review some standard facts about
conjugacy classes in compact Lie groups and 
extend them to the case of twisted conjugacy classes. A rather algebraic
approach relating twisted conjugacy classes to irreducible highest weight
representations of twisted affine Lie algebras via boundary conformal field
theory has been developed in \cite{BFS}, \cite{FS}.

The contents of this note is as follows: In section \ref{conj}, we
review some basic facts about twisted 
conjugation. In section \ref{conjcl}, we 
explicitly describe the set of twisted conjugacy classes of
a compact simply connected Lie group $G$ as
a convex polytope in an Euclidean vector space whose faces of maximal
dimension are in one to one correspondence with the vertices of a certain
twisted affine Dynkin diagram.
In fact, this polytope can naturally be 
identified with the fundamental domain of the twisted affine Weyl group
corresponding to $G$ and an automorphism $\tau$ of $G$. This
generalizes the classical situation of 
ordinary conjugacy classes which are parametrized by the fundamental alcove of
the affine Weyl group corresponding to $G$. In section \ref{stab}, we show how
to calculate the stabilizers of the twisted conjugacy classes. Since we
chose $G$ to be simply connected, the stabilizers of twisted 
conjugacy classes are connected subgroups of $G$. The Dynkin
diagram of the Lie algebras of the stabilizers turn out to be exactly 
the sub-diagrams
of the twisted affine Dynkin diagram described in section
\ref{conjcl}. Again, this is well known in the
untwisted case. 
The fundamental groups of the stabilizers can be easily calculated from
the root data so that we have a complete description of the stabilizers.

\medskip

Finally, section \ref{int} contains the main results of this paper. 
Suppose, the group $G$ is simple and and simply connected.
Let $\eta$ denote the Cartan 3-form on $G$. That is,  $\eta$ is a left
invariant 3-form on $G$ which 
generates $H^3(G,\Z)$. Let $\mathcal C$ denote an untwisted or twisted
conjugacy class in $G$, and let $\iota:\mathcal C\to G$ denote the
inclusion. It is known that the 3-form $\iota^*\eta$ on $\mathcal C$
is exact. Hence we can choose
a 2-form $\varpi$  on $\mathcal C$ such that $d\varpi=\iota^*\eta$.
Finally,
fix some $a\in\R^*$. We call 
an untwisted or twisted conjugacy classes $\mathcal C \subset G$ integral at
level $a$, if the integral
\begin{equation}\label{cond0}
  a\int_N \eta - a\int_{\partial N}\varpi
\end{equation} 
takes values in $\Z$ for all 3-cycles $N\subset G$ with 
$\partial N\subset\mathcal C$.
In other words, we call a conjugacy class $\mathcal C$ integral at
level $a$, if the relative 3-cocycle $a(\eta,\varpi)$ defines an
integral
element of the relative cohomology group $H^3(G,\mathcal C)$

It is not hard to see that
integral conjugacy classes at level $a$ can only exist if $a$ itself is an
integer.  
Furthermore, 
there is a well known correspondence between (twisted) conjugacy classes
in $G$ and coadjoint orbits of the corresponding (twisted) loop group $L(G)$
(\cite{F}, \cite{W}).
Let $\Omega G$ denote the space
of based loops in $G$. Using a 
modified version of the transgression homomorphism
$\sigma: H^3(G)\to H^2(\Omega G)$ we can pull back the relative
3-cocycle $a(\eta,\varpi)$ to a 2-cocycle on the coadjoint orbit of
$L(G)$ which corresponds to the twisted conjugacy
class $\mathcal C$. 
It turns out that this 2-cocycle is in the same cohomology class
as the standard symplectic structure on the coadjoint orbit. This
allows to
translate the integrality condition
(\ref{cond0}) for a
conjugacy class $\mathcal C$ to an integrality condition for the natural
symplectic form on the corresponding
coadjoint orbit of $L(G)$.
It is known 
\cite{F}, \cite{PS}, \cite{W}, that integral coadjoint orbits of a
(twisted) loop group correspond naturally to integrable irreducible
highest weight representations of the corresponding affine
Lie algebra. So we get a natural one-to-one correspondence between
integral conjugacy classes at non-negative level
and  integrable irreducible highest weight representations.


\section{Twisted conjugation}\label{conj}
Let $G$ be a Lie group and let $\tau$ be an automorphism of $G$.
The $\tau$--twisted conjugacy classes of $G$ are the orbits of the following
action of the group $G$ on itself: 
$$
  G\times G\to G,
$$ 
$$
  (g,h)\mapsto gh\tau(g^{-1})\,.
$$ 
Let $\tau'$ be another
automorphism of $G$ which differs from $\tau$ by an inner automorphism. That
is, $\tau'(g)=\tau(ugu^{-1})$ for some $u\in G$.  
Then the map $h\mapsto uh$
maps $\tau'$--twisted conjugacy classes in $G$ to $\tau$--twisted conjugacy
classes and induces an diffeomorphism between the corresponding twisted 
conjugacy classes. From this point of view, it is
enough to consider automorphisms $\tau$ up to inner automorphisms of $G$. 

Now let us suppose that $G$ is compact and semi-simple. Let $Aut(G)$ denote
the group of automorphisms of $G$, and let $Int(G)$ denote the subgroup of
inner automorphisms. Then $Aut(G)/Int(G)$ is a finite group. After fixing
a maximal
torus $T\subset G$ and a basis $\Pi$ of the root system $\Delta$ 
of $G$ with respect to $T$, the group $Aut(G)/Int(G)$  can be
identified with a subgroup of the
group of graph--automorphisms of the Dynkin diagram of the 
Lie algebra of $G$. In fact, if we
assume $G$ to be simply connected, then the group $Aut(G)/Int(G)$ is
isomorphic to the group of graph--automorphisms of the Dynkin diagram. This
observation allows us to find a nice representative in each connected component
of $Aut(G)$. Namely, after fixing a maximal torus $T\subset G$ and a basis
$\Pi$ of the root system $\Delta$ of $G$, any graph--automorphisms of the
Dynkin diagram corresponds to a permutation of the set $\Pi$. Furthermore, it
can be lifted to an automorphism of $G$ leaving $T$ invariant. From now on, we
will only consider automorphisms $\tau$ of $G$ which come from this
construction. 

One can view twisted conjugacy classes as $G$-orbits 
in the non-connected Lie group $\widetilde G=G\rtimes Aut(G)/Int(G)$, 
where $G$ acts on $\widetilde G$ by conjugation. 
From this point of
view, the twisted conjugacy classes have been studied in \cite{dS} 
(see also \cite{M} and  \cite{W}).


\section{The space of conjugacy classes}\label{conjcl}
Our first goal is to describe the set of $\tau$--twisted conjugacy
classes in $G$. 
Choosing $\tau$ to be the identity yields a description of the set of
ordinary conjugacy classes in $G$.
Since $\tau$ leaves the maximal torus $T\subset G$ invariant, we can consider
the subgroup $T^\tau\subset T$ of $\tau$--invariants. This group will in
general not be connected. Let us denote by $T^\tau_0$  the connected
component of $T^\tau$ containing the identity.
It is a  
fact that every $\tau$--twisted conjugacy class in $G$
intersects $T^\tau_0$ in at least one point (see e.g. \cite{BtD},
Proposition 4.3). 

So in order to describe the set of $\tau$--twisted conjugacy classes
in $G$, 
it
remains to check, which elements of  $T^\tau_0$ are twisted 
conjugate under $G$.
To this end, let us introduce the twisted Weyl group 
$$
  W(G,T,\tau)=N^\tau_{G}(T^\tau_0)/T^\tau_0\,.
$$
Here, $N^\tau_{G}(T^\tau_0)=\{g\in G~|~gT^\tau_0\tau(g)^{-1}=T^\tau_0\}$
denotes the normalizer of $T^\tau_0$ with respect to $\tau$-twisted
conjugation. 
It is a general fact  that $W(G,T,\tau)$ 
is a finite group (see e.g.\cite{BtD}). 
The twisted Weyl group $W(G,T,\tau)$ can be seen as a
generalization of  the ordinary Weyl group
$N(G,T)=N_G(T)/T$ of $G$ with respect to the maximal torus $T$.
One can show that two elements of $T^\tau_0$ are twisted conjugate 
under $G$ if
and only if they are conjugate under $W(G,T,\tau)$.
So the space of $\tau$-twisted conjugacy classes in 
$G$ can be identified with the
quotient $T^\tau_0/W(G,T,\tau)$. We shall now describe 
this quotient in
more detail.

Let $W=N_G(T)/T$ denote the usual Weyl group of $G$. The action of 
$\tau$ on the 
torus $T$ induces an action of $\tau$ on the Weyl group $W(G,T)$. 
Let $W^\tau$ denote the
subgroup of $W(G,T)$ which consists of elements commuting with $\tau$.
We also introduce the finite group $(T/T^\tau_0)^\tau$. Then one can show 
the following isomorphism 
(\cite{M}, \cite{W}):
$$
  W(G,T,\tau)\cong W(G,T)^\tau\ltimes (T/T^\tau_0)^\tau\,.
$$
The group $W^\tau$ is the Weyl group of the identity component of the
fixed point group $G^\tau$. Since  $T^\tau_0$ is a maximal torus of this group
and every element of $W^\tau$ commutes with $\tau$, the group $W^\tau$ acts on
$T^\tau_0$ by its Weyl group action on $T^\tau_0$. Now let us study the action
of $(T/T^\tau_0)^\tau$ on $T^\tau_0$. Take some element $\bar t\in
(T/T^\tau_0)^\tau$ and fix a pre-image $t$ of $\bar t$ under the projection
$T\to T/T^\tau_0$.   The condition that $\bar t$ is invariant under $\tau$
translates to the equation $\tau (t)=ts$ for some $s\in
T^\tau_0$. Hence $t t_0\tau (t^{-1})=t_0s$ for all $t_0\in T^\tau$.  Thus,
$T/T^\tau_0$ acts on $T^\tau_0$ by translations. 

The group $T/T^\tau_0$ and its action on $T^\tau_0$ can be described more
explicitly: Let $\h$ denote the Lie algebra of $T$. The Killing form on $G$
induces a $W(G,T)$--invariant inner product on $\h$ which is also invariant
under 
the action of $\tau$ on $\h$. Let $\h^{\tau}$ denote the $\tau$--invariant
part of $\h$, i.e. the Lie algebra of $T^\tau_0$ and let $\pi:\h\to\h^{\tau}$
denote the orthogonal projection with respect to the Killing form on
$\h$. Finally, let $\exp:\h\to T$ denote the exponential map. Its kernel is a
lattice in $\h$ and we can identify $T$ with $\h/ker(\exp)$. Since $\tau$ acts
on $T$, its induced action on $\h$ leaves the lattice $ker(\exp)$ invariant so
that we have $T^\tau_0=\h^{\tau}/ker(\exp)^\tau$. Finally, one checks
directly that $(T/T^\tau_0)^\tau\cong \pi(ker(\exp))/ker(\exp)^\tau$ and that
the translation action of $(T/T^\tau_0)^\tau$ on $T^\tau_0$ comes from the
translation action of $\pi(ker(\exp))$ on $\h^\tau$.

The observations above allow us to identify  the space of twisted conjugacy
classes with the set $\h^\tau/W^\tau\ltimes\pi(ker(\exp))$. 
This set can be identified with
a polytope in $\h^\tau$ as follows. First, let $K\subset\h$ denote the
fundamental Weyl chamber 
$K=\{h\in\h~|~\alpha(h)>0\text{ for all }\alpha\in \Pi\}$.
Then $K\cap\h^\tau\neq\emptyset$ and we have
$K\cap\h^\tau=\{h\in\h^\tau~|~\alpha|_{\h^\tau}(h)>0\text{ for all }
\alpha\in \Pi\}$.
The closure 
$\bar K=\{h\in\h~|~\alpha(h)\geq0\text{ for all }\alpha\in \Pi\}$ of
$K$ is a fundamental domain for the action of the 
Weyl group $W$ of $G$ on $\h$. Similarly, the closure $\bar K\cap\h^\tau$ is a
fundamental domain for the action of $W^\tau$ on $\h^\tau$. 

It remains to find a fundamental domain for the action of the group
$W^\tau\ltimes\pi(ker(\exp))$ on $\h^{\tau}$. It is well known that if the
group $G$ is simply connected and for
$\tau=id$, the lattice $ker(\exp)$ is the dual root lattice of $G$. Hence, the
group $W\ltimes ker(\exp)$ is the affine Weyl group of $\g$, or equivalently
the Weyl group of the affine Kac-Moody algebra corresponding to $\g$. 
A fundamental domain for the action of $W\ltimes ker(\exp)$ is
given by the fundamental alcove, i.e. by the set
$$
\mathfrak a=\{h\in\h~|~\alpha(h)\geq0\text{ for all }\alpha\in\Pi
\text{ and }\theta(h)\leq 1\}\,.
$$ 
Here,
$\theta\in\Delta$ denotes the highest root of $\Delta$ with respect to $\Pi$. 
If $G$ is not simply connected, 
we have to divide $\h$ by $W\ltimes \Lambda$, where
$\Lambda\subset\h$ denotes the lattice of all smooth 
homomorphisms $S^1\to T$.

We now want a similar description for a fundamental
domain for the action of $W^\tau\ltimes\pi(ker(\exp))$ on $\h^\tau$. 
First, we have to make a general observation. Denote by $\Delta^\tau$ the set 
$\Delta^\tau=\{\alpha|_{\h^\tau}~|~\alpha\in\Delta\}$. This is a subset in
$(\h^\tau)^*$, and if $\Delta$ is not of type $\text{A}_{2n}$, 
the set $\Delta^\tau$
is a root system. If $\Delta$ is of type $A_{2n}$ with $n>1$, 
then $\Delta^\tau$ is a
non-reduced root system of type $BC_n$, i.e. it is built out of the root
system $B_n$ and $C_n$ such that the long roots of $B_n$ are the short roots
of $C_n$. In the case $\Delta=A_2$, the set $\Delta^\tau$ consists of the
union of two root systems of type $A_1$ such that each element in one copy of
$A_1$ is two times an element of the other copy of $A_1$. 

We use the observation (\cite{W}) 
that the group $W^\tau\ltimes \pi(ker(\exp))$ is the Weyl
group of the twisted affine Lie algebra corresponding to $\g$ and the
automorphism $\tau$. A fundamental domain for the action of this group on
$\h^\tau$ can be described as follows:
If $\Delta$ is not of type $A_{2n}$, let $\theta_\tau$ 
denote the highest short root of $\Delta^\tau$ with respect to the basis
$\{\alpha|_{\h^\tau}~|~\alpha\in\Pi\}$. 
If $\Delta$ is of type $A_{2n}$ with $n>1$, 
let $\theta_\tau$ denote the highest short
root of the subsystem $B_n$ of $BC_n$ multiplied by two
(i.e. $\theta_\tau$ is a long root of
the system $C_n$). If $\Delta$ is of type $A_2$, let
$\theta_\tau$ denote the unique positive long root of $\Delta^\tau$. 

Now, in all cases, a fundamental domain for the
action of $W^\tau\ltimes\pi(ker(\exp))$ on $\h^\tau$ is given by the set 
$$
  \mathfrak a_\tau=\{h\in\h^\tau~|~\alpha|_{\h^\tau}(h)\geq0
  \text{ for }
  \alpha\in\Pi\text{ and }\theta_\tau(h)\leq\frac{1}{ord(\tau)}\}\,.
$$
A proof of this fact can be found e.g. in \cite{K}, chapter 6.
(Note that we have used a different  normalization of the invariant
bilinear form on the twisted affine Lie algebra than the one used in
\cite{K}.)  
Again, if $G$ is not simply connected, we have to divide $\h^\tau$ by 
the action of  $W^\tau\ltimes\pi(\Lambda)$.

\medskip

\begin{Example}
Let $G=SU(n)$ be the special unitary group. 
A maximal torus $T\subset SU(n)$ is given by the set of diagonal
matrices. The exponential map 
$$
  \exp:(x_1,\ldots x_n)\mapsto 
    diag(e^{2\pi i x_1},\ldots, e^{2\pi i x_n})\,,
$$
identifies the Lie algebra of $T$
with the set $\{(x_1,\ldots x_n)\in\R^n~|~\sum x_j=0\}$.
One can check directly, that the set of conjugacy classes in
$SU(n)$ is parametrized by the fundamental alcove
$$
  \mathfrak a=\{(x_1,\ldots,x_n)\in\R^n~|~x_1\geq\ldots\geq 
  x_n,~ \sum_ix_i=0,\text{ and }x_1-x_n\leq1\}\,.
$$
%
\end{Example}

\begin{Example}
For $n\geq3$, The group $SU(n)$ admits a non-trivial outer 
automorphism $\tau$
which is defined as follows.
Let $J_n$ denote the matrix
$J_n=\text{antidiag(1,\ldots,1)}$ if $n$ is odd and
$J_n=\text{antidiag(1,\ldots,1,-1,\ldots,-1)}$ if $n$ is even.
Then we can set $\tau(A)=J_n \bar A J_n^{-1}$ for $A\in SU(n)$. 
The automorphism $\tau$ acts on
$\h$ via
$\tau:diag(x_1,\ldots,x_n)\mapsto diag(-x_n,\ldots,-x_1)$. Hence
$\h^{\tau}=\{ diag(x_1,\ldots,x_n)~|~x_i=-x_{n+1-i}\text{ for }1\leq i\leq
n\}$.  
%
If $n=2m$ is even, the polytope $\mathfrak a_\tau$ 
which parametrizes the set of
$\tau$-twisted conjugacy classes in $SU(n)$ is given by
$$
  \mathfrak a_\tau = 
   \{(x_1,\ldots,x_m)\in\R^m~|~x_1\geq\ldots\geq 
x_m\geq0,\text{ and }x_1+x_m\leq\frac12\}\,.
$$
%
%
If $n=2m+1$ is odd and $m>1$, the polytope
$\mathfrak a_\tau$ is given by 
$$
  \mathfrak a_\tau = 
   \{(x_1,\ldots,x_m)\in\R^m~|~x_1\geq\ldots\geq 
  x_m\geq0,\text{ and }x_1\leq\frac14\}\,.
$$
Finally, in the case $SU(3)$, we have $\h^\tau=\{(x,0,-x)~|~x\in\R\}$, 
and
$$
 \mathfrak a_\tau=\{x\in\R~|~0\leq x\leq\frac14\}\,.
$$
\end{Example}


\section{Stabilizers}\label{stab}
Throughout this section, let $G$ be simply connected. Given an element $h\in
G$, the conjugacy class containing $h$ is isomorphic to $G/Stab_G(h)$, where
$Stab_G(h)=\{g\in G~|~gh\tau(g^{-1})=h\}$ denotes the stabilizer of
$h$ in 
$G$. The aim of this section is to give an explicit description of the
stabilizers. We first describe the Lie algebras of the stabilizer. 
Using this description along with some general facts of the theory
of compact Lie groups, we can calculate the fundamental groups of the
stabilizer. The fact that $G$ is a simply connected compact Lie group
implies that all
stabilizers are connected. So we have a complete description of the
twisted conjugacy classes in $G$.

As in section \ref{conj}, let $\tau$ denote an automorphism of $G$ which
leaves a maximal torus $T\subset G$ invariant and induces an automorphism of
the corresponding Dynkin diagram of $G$. Let $\mathfrak a\subset\h$ denote the
polytope parametrizing the set of conjugacy classes in $G$, and let $\mathfrak
a_\tau$ denote the polytope in $\h^\tau$ parametrizing the set of
$\tau$--twisted conjugacy classes in $G$.

Let us set $\widetilde\Pi=\Pi\cup\{-\theta\}$, where, as before, $\Pi$ denotes
a basis of $\Delta$. Then $\widetilde\Pi$
are the vertices of the extended Dynkin diagram of $\Delta$ or equivalently,
the Dynkin diagram of the affine Lie algebra corresponding to $G$. 
Similarly, the
set $\widetilde \Pi^\tau=\Pi^\tau\cup\{-\theta_\tau\}$ labels the vertices of
the Dynkin diagram of the twisted affine Lie algebra corresponding to $G$ and
the automorphism  $\tau$ of $G$.

\medskip

It is implicitly contained in the classical literature on compact Lie
groups (see e.g. \cite{A}), 
that
the stabilizers of elements 
of $G$ under ordinary conjugation have Dynkin diagrams which 
can be obtained by deleting vertices from the 
Dynkin diagram corresponding to $\widetilde\Pi$. 
This fact is a special case of the following proposition in the case $\tau=id$.
\footnote{The method described
in \cite{IS} does not give the complete set of the possible stabilizers.
E.g. the group B$_4$ appears as a stabilizer in F$_4$, but cannot be
obtained by deleting vertices of the unextended
  Dynkin
  diagram of F$_4$.}

\begin{Prop}\label{propstab}
  Let $G$ act on itself by $\tau$--twisted conjugation and let
  $H\in\mathfrak a_\tau$. Then the Dynkin diagram of the Lie algebra of
  $Stab(\exp(H))$ is the sub-diagram of the 
  Dynkin diagram corresponding to
  $\widetilde\Pi^\tau$ 
  which is obtained by deleting all $\alpha$ from the finite Dynkin diagram
  $\Pi^\tau$ for which
  $\alpha(H)\not\in\Z$ and deleting the vertex corresponding to $\theta_\tau$
  whenever $\theta_\tau(H)\not \in\frac{1}{ord(\tau)}+\Z$.

  The Lie algebra of $Stab(\exp(H))$ is given by
  the sum 
  of $\h^\tau$ and the sub-algebra of $\g$ corresponding to the diagram
  described above.   
\end{Prop}

\begin{proof}
Let us start with the case that $\g$ is not of type $A_{2m}$. 
The group $Stab_G(exp(H))$ can be written as
$Stab_G(exp(H))=\{g\in G~|~\exp(H)\tau(g)\exp(-H)=g\}$. 
Therefore its
Lie algebra is given by 
$$
  Lie(Stab_G(\exp(H)))=
  \{X\in\g~|~Ad(\exp(H))\circ\tau(X)=X\}\,.
$$
We can decompose the Lie algebra $\g_\C$ into the eigenspaces of $\tau$:  
$$ 
  \g_\C=\bigoplus_{k=0}^{\ord(\tau)-1}\g_k\,,
$$
where $\g_k$ denotes the $e^{2\pi i \frac{k}{ord(\tau)}}$--eigenspace of the
action of $\tau$ on $\g_\C$. It is known that $\g_0$ is a semi-simple Lie
algebra with root system $\Delta^\tau$ and that $\g_k$ with $k\neq0$ are
representations of $\g_0$ whose highest weight is given by $\theta_\tau$, the
highest short root of $\Delta^\tau$.  The Cartan sub-algebra of $\g_0$ is given
by 
$\h_\C^\tau$. Decomposing $\g_k$ into weight spaces with respect to $\h^\tau$,
we can write
$$
  \g_\C=\h^\tau\oplus\bigoplus_{\alpha\in\Delta^\tau}\g_\alpha
   \oplus\bigoplus_{k=1}^{ord(\tau)-1}\bigoplus_{\lambda\in P_k}\g_\lambda\,,
$$
where $P_k\subset(\h^\tau)^*$ 
denotes the set of weights of $\g_k$ as a representation of $\g_0$

Let us view $X$ as an element of $\g_\C$ and write
$$
  X=H_0+\sum_{\alpha\in\Delta^\tau}X_\alpha
   +\sum_{k=1}^{ord(\tau)-1}\sum_{\lambda\in P_k}X_\lambda\,,
$$
we see that we get 
$$
  Ad(\exp(H))\circ\tau(X)=H_0+\sum_{\alpha\in\Delta^\tau}e^{2\pi i\alpha(H)}
  X_\alpha
   +\sum_{k=1}^{ord(\tau)-1}\sum_{\lambda\in P_k}
   e^{2\pi i \frac{k}{ord(\tau)}}e^{2\pi i\lambda(H)}X_\lambda\,.
$$
This implies that the Lie algebra of $Stab(\exp(H))$ is given as the sum of
$\h^\tau_\C$ with those $\g_\alpha$ such that $\alpha(H)\in\Z$ and those
$\g_\lambda$ with $\lambda\in P_k$ such that
$\lambda(H)\in\frac{-k}{ord(\tau)}+\Z$.
Since we have chosen $H\in\mathfrak a_\tau$, and since $\theta_\tau$ is
the highest weight of the $\g_k$ we have $0\leq\alpha(H)\leq1$ and
$0\leq\lambda(H)\leq\frac{1}{ord(\tau)}$ for all $\alpha\in\Delta_\tau$ and
$\lambda\in P_k$. 
Furthermore,  $|\lambda(H)|=\frac{1}{ord(\tau)}$ can only be obtained for
 $\lambda=\pm\theta_\tau$ in which case $H$ has to lie on the boundary
 of  $\mathfrak a_\tau$. 
If $H$ lies in the interior of $\mathfrak a_\tau$, we have $\alpha(H)\not\in\Z$
 for all $\alpha\in\Delta^\tau$ and
 $\lambda(H)\not\in\frac{-k}{ord(\tau)}+\Z$ for all
 $\lambda\in P_k$. 
If $H$ lies on the boundary of $\mathfrak a_\tau$, 
the stabilizer of $H$ is generated by the elements of
$\widetilde\Pi^\tau=\Pi^\tau\cup\{\theta_\tau\}$ 
for which either $\alpha(H)\in\Z$ if
$\alpha\in\Pi^\tau$ or $\theta_\tau(H)\in\frac{-k}{ord(\tau)}+\Z$.

The case that $\g$ is of type $A_{2m}$ can be obtained by similar arguments,
one just has to be more careful with the root system $\Delta^\tau$ which is
non-reduced.  
\end{proof}

\begin{Remark}
There is a second approach to Proposition \ref{propstab} 
which works for semi-simple conjugacy classes in algebraic
groups. This approach uses the classification of finite order
automorphisms of simple Lie algebras and has been worked out in the 
proof of Theorem 3.2 in \cite{M}.
\end{Remark}

Finally, we have to study the topology of the stabilizers. Since we
assumed $G$ to be simply connected, all stabilizers are
connected. So in order to describe them explicitly, we only have to determine
their fundamental groups. This can be done using some standard facts from the
theory of compact Lie groups (see e.g. \cite{BtD}):
Let $\Delta_H\subset(\h^\tau)^*$ denote the root system of the
group 
$Stab(\exp(H))$, let $Q_H\subset(\h^\tau)^*$ denote the lattice generated by
$\Delta_H$.

We can use the normalized Killing form on $\h$ to identify $(\h^\tau)^*$
with $\h^\tau$. Let $Q_H^\vee$ be the image of $Q_H$ under this
embedding.
Then
$$
  \pi_1(Stab(\exp(H)))\cong ker(\exp)^\tau/Q_H^\vee\,.
$$

\begin{Example}
In figure \ref{C2A3} 
we indicate how the constructions
described in this section apply to the spaces of
conjugacy classes of the simply connected compact Lie groups of type
$C_2$
compared to the space of twisted conjugacy classes of the simply
connected Lie groups of type $A_3$. 
The picture shows the fundamental domains of the respective untwisted and 
twisted affine
Weyl groups. The faces of the fundamental domain are labeled with the
stabilizers of the corresponding conjugacy classes.
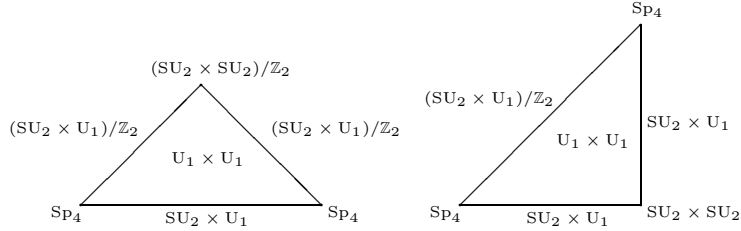
\begin{figure}[h]
\begin{center}
\setlength{\unitlength}{.8mm}
\begin{picture}(120,28)
\thinlines 
\put(2,-2){\tiny$\mbox{Sp}_4$}\put(48,-2){\tiny$\mbox{Sp}_4$}
\put(18.5,22){\tiny$(\mbox{SU}_2\times\mbox{SU}_2)/\mathbb{Z}_2$}
\put(38.5,12){\tiny$(\mbox{SU}_2\times\mbox{U}_1)/\mathbb{Z}_2$}
\put(-5,12){\tiny$(\mbox{SU}_2\times\mbox{U}_1)/\mathbb{Z}_2$}
\put(21,-3){\tiny$\mbox{SU}_2\times\mbox{U}_1$}
\put(22,7){\tiny$\mbox{U}_1\times\mbox{U}_1$}
\put(7,0){\circle*{.05}}\put(47,0){\circle*{.05}}\put(27,20){\circle*{.05}}
\put(7,0){\line(1,0){40}} \put(7,0){\line(1,1){20}} 
\put(27,20){\line(1,-1){20}} 
\put(65,-2){\tiny$\mbox{Sp}_4$}\put(101,-2){\tiny$\mbox{SU}_2\times\mbox{SU}_2$}
\put(98.5,32){\tiny$\mbox{Sp}_4$}
\put(101,13){\tiny$\mbox{SU}_2\times\mbox{U}_1$}
\put(64,17){\tiny$(\mbox{SU}_2\times\mbox{U}_1)/\mathbb{Z}_2$}
\put(81,-3){\tiny$\mbox{SU}_2\times\mbox{U}_1$}
\put(86,10){\tiny$\mbox{U}_1\times\mbox{U}_1$}
\put(70,0){\circle*{.05}}\put(100,0){\circle*{.05}}
\put(100,30){\circle*{.05}}
\put(70,0){\line(1,0){30}} \put(70,0){\line(1,1){30}} 
\put(100,0){\line(0,1){30}} 
\end{picture}
\caption{Conjugacy classes of $Sp(4)$  (left) versus 
        $\tau$-twisted conjugacy classes of $SU(4)$ (right) where $\tau$ 
        is a a non-trivial outer automorphism of $SU(4)$}
\label{C2A3}
\end{center}
\end{figure}
\end{Example}


\section{Integral conjugacy classes}\label{int}
%


\subsection{The integrality condition}
We start this section with our main definition: Suppose that $G$ is
simply connected and simple and
let $\eta$ denote the left-invariant 3-form
on $G$ whose value at the identity is given by 
$\eta(X_1,X_2,X_3)=\frac{1}{8\pi^2}\langle X_1,[X_2,X_3]\rangle$. 
Thus, $\eta$ is a generator of $H^3(G,\Z)$. 
Let $\tau$ be an
automorphism of $G$,  let $\mathcal C$ be a $\tau$-twisted
conjugacy class on $G$ and denote by $\iota:\mathcal C\to G$ the
embedding of $\mathcal C$ into $G$. 
On $\mathcal C$, there is a canonical 2-form $\varpi$ on $\mathcal C$, 
which is defined as follows. Given an element $A\in\g$, denote by 
$A_{\mathcal C}$ the generating vector field of the action of $G$ on
$\mathcal C$. Then we set 
$$
  \varpi_g(A_{\mathcal C},B_{\mathcal C}) =
  \frac{1}{8\pi^2}\left(\langle \Ad_g\circ\tau(A),B\rangle
    - \langle\Ad_g\circ\tau(B),A\rangle\right)\,.
$$
It is known (see e.g. \cite{AMM} Proposition 2.1, \cite{AMW}) that
$$
  d\varpi=\iota^*\eta\,.
$$
\begin{Definition}
  Fix some $a\in\R^*$. We call a twisted or untwisted conjugacy class
  $\mathcal C\subset G$ integral at level $a$  if the integral
  \begin{equation}\label{cond}
    a\int_N \eta-a\int_{\partial N}\varpi
  \end{equation} 
takes values in $\Z$ for all 3-cycles $N\subset G$ with 
$\partial N\subset\mathcal C$.
\end{Definition}
In other words, we call a conjugacy class $\mathcal C$ integral at
level $a$ if the relative cocycle 
$a(\eta,\varpi)\in H^3(G,\mathcal C)$ is integral.

The goal of this section is to classify integral conjugacy classes in terms
of integral coadjoint orbits of certain centrally extended loop groups, and
thus to relate them to integrable highest weight representations of certain
affine Lie algebras.


\subsection{Coadjoint orbits of loop groups}\label{coad}
Let $G$ be as before. We denote by $L(G)=C^\infty(S^1,G)$ the loop group of
$G$.\footnote{For technical reasons it is often more convenient
  to consider the Banach Lie group of maps $S^1\to G$ 
 of some fixed Sobolev class $s > 1/2$. We
 will ignore such subtleties throughout this note and stick with the
 more intuitive group of smooth loops}
 There is a universal central extension $\widehat{L(G)}$
of the group $L(G)$ which, as a topological space, is
a non-trivial $S^1$-bundle over $L(G)$. The underlying vector space of
the Lie algebra of $\widehat{L(G)}$
is given by $\widehat{L(\g)}=L(\g)\oplus\R$, where $L(\g)=C^\infty(S^1,\g)$
denotes the loop algebra of $\g$. The smooth part of the dual of
$\widehat{L(\g)}$ can be identified with $L(\g)\oplus\R$ via the pairing
\begin{equation}\label{pair}
  \langle (X,a),(Y,b)\rangle = \frac{1}{2\pi}\int_{S^1}\langle
  X(\theta),Y(\theta)\rangle d\theta + ab\,.
\end{equation}
Here, the bilinear form $\langle~,~\rangle$ on the right hand denotes the
Killing form on the Lie algebra $\g$ normalized so that the long roots of
$\g_\C$ have square length $2$.

The coadjoint action of the loop group $L(G)$ on $L(\g)\oplus\R$ is
given by \cite{PS}, \cite{F} 
$$
  \Ad_g^*:(X,a)\mapsto (gXg^{-1}-ag'g^{-1},a)\,,
$$
where $g$ is an element of $L(G)$ and $g'=\frac{d}{d\theta}g$ denotes the
derivative of $g$ with respect to $\theta$.
Let us denote by $\mathcal O_{(X,a)}\subset L(\g)\oplus\R$ 
the $L(G)$--orbit through
$(X,a)$. 

For $a\neq 0$, we can solve the
differential equation
$$  
z' = -\frac1a Xz\,.
$$
Let us denote by $z_{(X,a)}$ the unique solution of this equation with initial
condition $z_{(X,a)}(0)=e$, where  $e$ is the identity element in
$G$. This is a path in the Lie group $G$ starting at
$e$. Now, since $X\in L(\g)$ is periodic, 
we get $z_{(X,a)}(\theta+2\pi)=z_{(X,a)}(\theta) z_{(X,a)}(2\pi)$.
Taking another element $(Y,a)$ in the coadjoint orbit
$\mathcal O_{(X,a)}$, we have $Y=gXg^{-1}-ag'g^{-1}$ for some
$g\in L(G)$. An easy calculation \cite{PS}, \cite{F}
shows that 
\begin{equation}\label{gz}
  z_{(Y,a)}(\theta)=g(\theta)z_{(X,a)}(\theta)g(0)^{-1}\,.
\end{equation}
Since $g$ is periodic, we see that $z_{(X,a)}(2\pi)$ and $z_{(Y,a)}(2\pi)$ 
lie in the same
conjugacy class in $G$. Furthermore, the stabilizer of $(X,a)$ in $L(G)$ is
isomorphic to the stabilizer of $z_{(X,a)}(2\pi)$ 
in $G$ so that we get 
$\mathcal O_{(X,a)}\cong L(G)/Stab_G(z_{(X,a)}(2\pi))$. 
The stabilizers have been explicitly
described in section \ref{stab}.

Now, let $\tau$ be an automorphism of $G$ of finite order $ord(\tau)=r$. 
The twisted loop group $L(G,\tau)$  
is defined as follows: $L(G,\tau)=\{g\in
C^\infty(\R,G)~|~g(\theta)=\tau(g(\theta+2\pi))\}$.
As in the untwisted case, there is a universal central extension
$\widehat{L(G,\tau)}$ of
$L(G,\tau)$ by the circle group $S^1$ which, as a topological space, is a
non-trivial $S^1$--bundle over $L(G,\tau)$. The Lie algebra of
$\widehat{L(G,\tau)}$ is  $\widehat{L(\g,\tau)}=L(\g,\tau)\oplus\R$, where
$L(\g,\tau)=\{X\in C^\infty(\R,\g)~|~X(\theta)=\tau(X(\theta+2\pi))\}$ 
denotes the
twisted loop algebra of $\g$. Again, the smooth part of the dual of
$\widehat{L(\g,\tau)}$ can be identified with $L(\g,\tau)\oplus\R$ via a
non-degenerate pairing 
$$
\langle(X,a),(Y,b)\rangle = \frac{1}{2\pi}\int_0^{2\pi}\langle
X(\theta),Y(\theta)\rangle d\theta +ab\,.
$$ 
The coadjoint action of $L(G,\tau)$ on
$L(\g,\tau)\oplus\R$ is again given by $g:(X,a)\mapsto
(gXg^{-1}-a g'g^{-1},a)$.
It has been observed in \cite{W} that the correspondence between coadjoint
orbits of $L(G)$ and conjugacy classes in $G$ extends to a correspondence
between coadjoint orbits of $L(G,\tau)$ and $\tau$--twisted conjugacy classes
in $G$ as follows: Fix $(X,a)\in L(\g,\tau)\oplus\R$ with $a\neq 0$ 
and, as before, solve the
differential equation $z'=-\frac1a Xz$ 
with initial condition $z_{(X,a)}(0)=e$. 
If $(Y,a)$ lies in the same coadjoint orbit as $(X,a)$, a similar calculation
as in the untwisted case shows that $z_{(X,a)}(2\pi)$ and 
$z_{(Y,a)}(2\pi)$ lie in
the same $\tau$--twisted conjugacy class. 
Furthermore, the stabilizer of $(X,a)$ in
$L(G,\tau)$ is isomorphic to the stabilizer of $z_{(X,a)}(2\pi)$ 
in $G$
(where we 
mean the stabilizer with respect to $\tau$--twisted conjugation in $G$).
So we have 
$\mathcal O_{(X,a)}\cong L(G,\tau)/Stab(z_{(X,a)}(2\pi))$.

On every coadjoint orbit of a Lie group there exists a natural symplectic
structure. In our situation, the corresponding 2-form $\omega$ on $\mathcal
O_{(X,a)}$ with $a\neq 0$ is given as
follows:  Fix some  $(Y,a)\in\mathcal O_{(X,a)}$ 
and let $A_1,~A_2$ be two tangent vectors at $(Y,a)$. We can view $A_i$ as
elements of $L(\g)$.
Then the 2-form $\omega$ at $(Y,a)$ evaluated at $A_1$ and $A_2$ is given by
$$
  \omega_{(Y,a)}(A_1,A_2) = \frac{1}{2\pi}\int_{S^1}\langle
  Y(\theta),[A_1(\theta),A_2(\theta)]\rangle d\theta +
  \frac{a}{2\pi}\int_{S^1}\langle A_1'(\theta), A_2(\theta)\rangle d\theta \,.
$$

For twisted loop groups the 2-form on the coadjoint orbits is given by the
same formula.


\subsection{Integral conjugacy classes}
As we have seen in the last section, we can identify each element 
$(Y,a)\in \mathcal
O_{(X,a)}$ with the unique solution of the differential equation
$z'=-\frac1a Yz$
with initial condition $z(0)=e$. 
This allows to define a map
$$
F:\R\times\mathcal O_{(X,a)}\to G
$$ 
via 
$$
  F:\left(\theta,(Y,a)\right)\mapsto z_{(Y,a)}(\theta)\,.
$$ 
As before, let $\eta$ denote the closed 
left-invariant 3-form on $G$ whose value at the
identity element of $G$ evaluated on three tangent vectors $X_1,X_2,X_3\in\g$
is given by 
$\eta(X_1,X_2,X_3)=\frac{1}{8\pi^2}\langle X_1,[X_2,X_3]\rangle$. 

Using the map $F$, we can pull back $\eta$ to a 3-form $F^*\eta$ on
$\R\times\mathcal O_{(X,a)}$. Its value at
$(\theta,(Y,a))\in\R\times\mathcal O_{(X,a)}$ evaluated on a triple
$(\delta \theta,\delta_1 (Y,a),\delta_2(Y,a))$ of tangent vectors at
$(\theta,(Y,a))$
is 
$$
  \frac{1}{8\pi^2}\langle z_{(Y,a)}'(\theta)z_{(Y,a)}(\theta)^{-1},
  [\xi_1(\theta),\xi_2(\theta)]\rangle\,,
$$
where we have identified the tangent vectors $\delta_i(Y,a)$ with vector
fields $\delta_iz_{(Y,a)}$ along $z_{(Y,a)}$ and have set 
$\xi_i(\theta) = \delta_iz_{(Y,a)}(\theta)z_{(Y,a)}(\theta)^{-1}$.
In particular, it follows from equation (\ref{gz}) that
in this identification, the value of the generating
vector field $B_{\mathcal O}$ corresponding to an element $B\in L(\g)$ at 
$z=z_{(Y,a)}$ is given by
$$
  B_{\mathcal O}(z)=(r_z)_*B-(l_z)_* B(0)\,,
$$
Where $r_g$ and $l_g$ denote right and left translation by an element 
$g\in G$.
In particular, we find 
$$
  \xi_B(\theta) = B(\theta) - Ad_{z(\theta)} B(0)\,.
$$

Now we can
integrate the form $F^*\eta$ over $[0,2\pi]$ to obtain a 2-form
$$
  \sigma(\eta)=\int_{0}^{2\pi}F^*\eta(\theta)d\theta
$$ 
on the coadjoint orbit $\mathcal O_{(X,a)}$. 

Let $\mathcal C$ denote the (twisted) conjugacy class corresponding to
the coadjoint orbit $\mathcal O_{(X,a)}$.
We can define a map $F_{2\pi}:\mathcal O_{(X,a)}\to\mathcal C$ via
$$
  F_{2\pi}(Y,a)= z_{(Y,a)}(2\pi)
$$
and use this map to pull back the 2-form $\varpi$ on $\mathcal C$ to
a 2-form $F_{2\pi}^*\varpi$ on $\mathcal O_{(X,a)}$. It is easy to see
that the two form 
$\sigma(\eta)-F_{2\pi}^*\varpi$ is closed. 

\begin{Prop}
  Let $a\neq 0$. Then the 2-form $a\sigma(\eta)-aF_{2\pi}^*\varpi$ 
  on the coadjoint orbit
  $\mathcal O_{(X,a)}$  
  is co-homologous to a multiple $\frac{1}{2\pi}\omega$ 
  of the Kirillov-Kostant
  symplectic form  on the coadjoint orbit $\mathcal
  O_{(X,a)}$. 
\end{Prop}

\begin{proof}
  We have to compare the two form  $2\pi a(\sigma(\eta)-F_{2\pi}^*\varpi)$
  with the the symplectic form $\omega$. 
  To this end, let us introduce a 1-form 
  $\beta$ on $\mathcal O_{(X,a)}$, which,
  at a point $(Y,a)\in \mathcal O_{(X,a)}$ evaluated at a tangent vector
  $\delta(Y,a)$ is given by
  $$
  \beta_{(Y,a)}(\delta(Y,a))=\frac{a}{4\pi}
  \int_{0}^{2\pi}\langle z'_{(Y,a)}(\theta)z_{(Y,a)}(\theta)^{-1},
    \xi(\theta)\rangle d\theta \,,
  $$
  with $\xi$ as before. Since we have
  $$\xi_i(\beta(\delta_j(Y,a))(Y,a)=\frac{a}{4\pi}\int_{0}^{2\pi}\langle
  \xi_i'(\theta),\xi_j(\theta)\rangle d\theta$$ 
  and 
  $$
    z_{(Y,a)}'(\theta)z_{(Y,a)}(\theta)^{-1}=-\frac1aY(\theta)\,,
  $$ 
  we find
  \begin{multline*}
    d\beta_{(Y,a)}(\xi_1,\xi_2) = \frac{1}{4\pi}\int_{0}^{2\pi}\langle
    Y(\theta),[\xi_1(\theta),\xi_2(\theta)]\rangle d\theta \\
    + 
    \frac{a}{2\pi}\int_{0}^{2\pi}\langle \xi_1'(\theta),
     \xi_2(\theta)\rangle d\theta - 
     \frac{a}{4\pi}\langle\xi_1(2\pi),\xi_2(2\pi)\rangle \,
  \end{multline*}
  and hence
  \begin{multline*}
    2\pi a\sigma(\eta)(\xi_1,\xi_2)+d\beta(\xi_1,\xi_2) = 
    \frac{1}{2\pi}\int_{0}^{2\pi}\langle
    Y(\theta),[\xi_1(\theta),\xi_2(\theta)]\rangle d\theta \\
    +
    \frac{a}{2\pi}\int_{0}^{2\pi}\langle \xi_1'(\theta), 
      \xi_2(\theta)\rangle d\theta - 
       \frac{a}{4\pi}\langle\xi_1(2\pi),\xi_2(2\pi)\rangle\,.
  \end{multline*}
  Since we can assume the $\xi_i$ to come from generating vectorfields of
  the $L(G)$-action, we can write 
  $\xi_i(\theta)=B_i(\theta)-Ad_{z(\theta)}B_i(0)$ for some $B_i\in
  L(\g)$ (respectively $B_i\in L(\g,\tau)$ in the twisted case) and
  $z=z_{(Y,a)}$.
  Inserting this into the equation above, a short calculation gives
  \begin{multline*}
    2\pi a\sigma(\eta)(\xi_1,\xi_2)+d\beta(\xi_1,\xi_2) = \\
      \frac{1}{2\pi}\int_{0}^{2\pi}\langle 
        Y(\theta),[B_1(\theta),B_2(\theta)]\rangle d\theta 
     +
     \frac{a}{2\pi}\int_{0}^{2\pi}\langle B_1'(\theta), 
     B_2(\theta)\rangle d\theta \\
     + 
     \frac{a}{4\pi}\langle Ad_{z(2\pi)}B_1(0),B_2(2\pi)\rangle
     -
    \frac{a}{4\pi}\langle Ad_{z(2\pi)}B_2(0),B_1(2\pi)\rangle\,.
  \end{multline*}
  Recall that $B_i(0)=\tau(B_i(2\pi))$, (with $\tau=id$ in the untwisted
  case) so that we have 
  $$
    \frac{a}{4\pi}\langle Ad_{z(2\pi)}B_1(0),B_2(2\pi)\rangle
     -
    \frac{a}{4\pi}\langle Ad_{z(2\pi)}B_2(0),B_1(2\pi)\rangle
  = 2\pi a F_{2\pi}^*\varpi(\xi_1,\xi_2)
  $$
  Hence we get 
  $$
    2\pi a\sigma(\eta) - 2\pi a  F_{2\pi}^*\varpi = \omega\,.
  $$
\end{proof}

Our next goal is to translate the integrality condition (\ref{cond})
for the pair $(\mathcal C,a)$ to an
integrality condition for the 2-form 
$a\sigma(\eta)-aF^*_{2\pi}\varpi$, and hence for the symplectic
form $\frac{1}{2\pi}\omega$ on the corresponding
coadjoint orbit $\mathcal O_{(X,a)}$. 
Given any closed 2-cycle $\widetilde N\in H_2(\mathcal O_{(X,a)})$, we
get a 3-cycle $N$ in $G$ by mapping 
$(Y,a)\in\widetilde N$ to
$\{z_{(Y,a)}(\theta)~|~0\leq\theta\leq2\pi\}$. By the results of the
last section, we have $\partial N\subset \mathcal C$ so that $N$ is
indeed a relative 3-cycle. Furthermore, by construction we
have
$$
  \int_{\widetilde N}\frac{1}{2\pi}\omega =
  \int_{\widetilde N}a(\sigma(\eta) - F_{2\pi}^*\varpi)= 
   a\int_N\eta - a\int_{\partial N} \varpi\,,
$$
so that integrality of $\frac{1}{2\pi}\omega$ 
is necessary for the integrality for $a\eta$. 

In the other direction, we need the following proposition.
\begin{Prop}\label{cycle}
For any relative 3-cycle $N$ in $G$ with $\partial N\subset\mathcal C$, 
there exists a (not necessarily unique) closed 
2-cycle $\widetilde N^*\in H_2(\mathcal O_{(X,a)})$ 
such that $\partial N^*=\partial N$. 
Here, the 3-cycle $N^*$ is obtained
from $\widetilde N^*$ by the construction described above.
\end{Prop}
We will postpone the proof of proposition \ref{cycle} to section
\ref{proof}.

\medskip 

Let us fix a 3-cycle $N$ in $G$ with $\partial N\subset\mathcal C$ 
and let $\widetilde N^*\in
H_2(\mathcal O_{(X,a)})$ be the 2-cycle from proposition \ref{cycle}.
As we have noted, the cycle $\widetilde N^*$ need not be
unique, but we can do with the following. 
Since  $\partial(N-N^*)=0$, we have $N-N^*\in H_3(G)$. 
We can write
$$
   a\int_N\eta-a\int_{\partial N}\varpi=
   a\int_{N-N^*}\eta + a\int_{N^*}\eta -a\int_{\partial N}\varpi
    = a\int_{N-N^*}\eta +
       a\int_{\widetilde N^*}
        \left(\sigma(\eta)-F_{2\pi}^*\varpi\right)\,.
$$

Now, let us take $N$ to be a generator of $H_3(G)\cong \Z$. Since 
$\partial N=\emptyset$, we find $N^*=0$. We have chosen $\eta$ 
to be the generator of
$H^3(G,\Z)\cong\Z$ so that the integrality condition (\ref{cond})
translates to 
$ 
  a\in\Z\,.
$
But this implies $a\int_{N-N^*}\eta\in\Z$ for all 3-cycles $N$ in $G$ with
$\partial N\subset\mathcal C$ so that integrality of
$a\sigma(\eta)-aF_{2\pi}^*\varpi$ implies integrality of $a(\eta,\varpi)$.
Putting everything together we get
\begin{Theorem}
  A pair $(\mathcal C,a)$ of an untwisted respectively twisted conjugacy 
  class $\mathcal C\subset G$ and an $a\neq0$
  satisfies the integrality
  condition (\ref{cond}) if and only if $a\in\Z\setminus\{0\}$ and the
  natural symplectic structure $\frac{1}{2\pi}\omega$ on 
  the corresponding coadjoint
  orbit $\mathcal O_{(X,a)}$ of the loop group $L(G)$, respectively the
  twisted loop group
  $L(G,\tau)$ is integral. 
\end{Theorem}

The integrality condition for the Kirillov-Kostant form 
$\frac{1}{2\pi}\omega$ on the
coadjoint orbit $\mathcal O_{(X,a)}$ can be translated back to an explicit
condition for the conjugacy class. Indeed, as we have seen in section
\ref{coad}, we can take $\frac 1aX$ to be a constant  
$\frac 1aX\in\mathfrak a\subset
\h\subset\g\subset L(\g)$ in the untwisted case and 
$\frac 1aX\in\mathfrak a_\tau\subset\h^\tau\subset\g^\tau
\subset L(\g,\tau)$ in the
twisted case. Now, the condition that $\frac{1}{2\pi}\omega$ 
is an integral 2-form
translates to the condition that $\alpha(X,a)\in\Z$ for all roots $\alpha$ of
the (twisted) affine Lie algebra $\widehat{L(\g_\C,\tau)}$. This gives again
the condition that $a$ must be an integer. Furthermore, let us
identify $\h\oplus\R$ with its dual via the non-degenerate pairing
from (\ref{pair}). Then, for fixed positive
$a\in\Z$, the $(X,a)$ satisfying the integrality condition are exactly the
highest weights of the irreducible highest weight representations of the
(twisted) affine Lie algebra $\widehat{L(\g,\tau)}$ at level $a$.

\begin{Remark}
In fact, one can associate to each irreducible highest weight representation a
coadjoint orbit of the corresponding (twisted) loop group \cite{F}, \cite{W}. 
Under this
correspondence, a highest weight representation with highest weight $(X,a)$
does not correspond to the coadjoint orbit passing through $(X,a)$ but rather
to the orbit through $(X+\rho_\tau,a+h_\tau^\vee)$, where $\rho_\tau\in\h^\tau$
denotes the projection of the half sum of all positive roots of $\g_\C$ to
$\h^\tau$, and $h^\vee_\tau$ denotes the dual Coxeter number of the twisted
affine Lie algebra $\widehat{L(\g,\tau)}$ (in the untwisted case, just take
$\tau=id$. See e.g. \cite{K} for more information on affine Lie algebras).
\end{Remark}


\subsection{Proof of Proposition \ref{cycle}}\label{proof}
Before we
start with the proof of Proposition \ref{cycle}, we need some preparations.
\begin{Lemma}\label{surj}
  For any $a\neq0$, the map $\mathcal O_{(X,a)}\to\mathcal C$ given by
  $(Y,a)\mapsto z_{(Y,a)}(2\pi)$ is surjective. 
\end{Lemma}
\begin{proof}
  First, note that the map $L(G,\tau)\to G$ which maps  $\gamma\in L(G,\tau)$
  to $\gamma(0)$ is surjective. Indeed, fix some $g\in G$ and choose $H\in\g$
  such that $\exp(H)=g$. We can decompose 
  $\g_\C=\bigoplus_{k=0}^{ord(\tau)-1}\g_k$, where $\g_k$
  denotes the $e^{2\pi i\frac{k}{ord(\tau)}}$--eigenspace of
  $\tau$. Viewing $H$ as an element of $\g_\C$,
  we can write $H=\sum_{k=0}^{ord(\tau)-1}H_k$. Let us set 
  $$
    Y_\C(\theta)=\sum_{k=0}^{ord(\tau)-1} H_k e^{\frac{ik}{ord(\tau)}\theta}\,
  $$
  The real part $Y$ of $Y_\C$ is an element of $L(\g,\tau)$ and we
  have $\exp(Y)(0)=g$.

  Finally, fix 
  some $g_0\in \mathcal C$. We can write $g_0=\tau(g)hg^{-1}$ for some
  $g\in G$ and $h\in G^\tau$. Fix $H\in \g^\tau$ with 
  $\exp(-\frac{2\pi}{a}H)=h$,
  and $\gamma\in L(G,\tau)$ with $\gamma(0)=g$. Set
  $Y_0(\theta)=\gamma(\theta)H\gamma(\theta)^{-1}-
   a\gamma'(\theta)\gamma^{-1}(\theta)$. Then we have 
  $z_{(Y_0,a)}(2\pi)=\gamma(2\pi)\exp(-\frac{2\pi}{a} H)\gamma(0)^{-1}=g_0$.
\end{proof}

\begin{proof}[Proof of Proposition \ref{cycle}]
%
%
Let us fix a triangulation $\{D_i\}_{i\in I}$ of $\partial N$ and a
point $x_i$ in the interior of each triangle. If
the triangulation is fine enough, each $D_i$ allows to choose 
a triangle $E_i$ in
the Lie group $G$ such that the identity element $e\in G$ is in the
interior of $E_i$ and we have $D_i=\{g x_i g^{-1}~|~g\in E_i\}$. Now,
using the fact that the exponential map $\exp:\g\to G$
is locally invertible at $e$,
we can refine the triangulation $\{D_i\}_{i\in I}$ to obtain triangles
$E_i'$ in $\g$ with $0\in E_i'$ such that $\exp$ restricted to an open
neighborhood of $E_i'$
invertible, and $\exp(E_i')=E_i$. By the construction in the first
part of the proof of 
Lemma \ref{surj}, the triangles $E_i'$ in $\g$ 
give triangles $E_i''\subset L(G,\tau)$ in the loop group.
Finally, using Lemma \ref{surj}, for
each $i\in I$ we
can choose $X_i\in\mathcal O_{(X,a)}$ such that 
$z_{(X_i,a)}(2\pi)=x_i$. Then the set 
$$
  D_i' = 
    \{\Ad^*_\gamma(X_i,a)~|~\gamma\in E_i''\}\subset \mathcal O_{(X,a)}
$$  
is a triangle in $\mathcal O_{(X,a)}$.

Now the idea is to take the union of all $D_i'$ with $i\in I$. 
But the triangles $D_i'$ 
might not fit together to form a closed 2-cycle in 
$\mathcal O_{(X,a)}$. Indeed, if an element
$g\in \partial N$ lies on the boundary of two triangles $D_1$ and $D_2$, 
the construction from
above associates two elements $(Y_1,a)$ and
$(Y_2,a)$ of $\mathcal O_{(X,a)}$ to
$g$ which might be different. 
But we know that both  $(Y_j,a)$ associated to such $g$ satisfy
$z_{(Y_j,a)}(2\pi)=g$. Set 
$$
  \mathcal F_g=\{(Y,a)\in\mathcal O_{(X,a)}~|~z_{(Y,a)}(2\pi)=g\}\,,
$$
and let us assume for the moment that 
$\pi_0(\mathcal F_g)=\pi_1(\mathcal F_g)=\{0\}$ for each $g\in G$. 

Then, since $\mathcal F_g$ is connected, we can join $(Y_1,a)$
and $(Y_2,a)$ by a path inside $\mathcal F_g$. Furthermore, we can
choose such a path for each point $g$ of the edge in such a way that
it depends continuously on the point $g$. 
Indeed,  
we can consider the union of all $\mathcal F_g$ with $g$ an element of a fixed
edge of the triangle $D_i$. 
This set is a subset of the coadjoint orbit and it is connected and
simply connected (being a fibration over a closed interval with fibers
isomorphic to $\mathcal F_g$).
This shows that that one can join the edges continuously (by contracting
the loop given by the edges together with the paths joining the endpoints
of the edges).

Doing
this for all edges of the triangulation $\{D_i\}_{i\in I}$, we can
thus ``join the edges'' of the $D_i'$. This procedure might still leave
``holes'' at the vertices of the triangulation. 
Let $g$ be a
vertex of the triangulation . 
The boundary of a ``hole at $g$'' is homeomorphic to an image
of $S^1$  inside $\mathcal F_g$. But since we have assumed $\mathcal F_g$ to
be simply connected, we can contract the boundary inside $\mathcal F_g$
and thereby ``fill the hole''. Repeating this process at each vertex of the
triangulation, we obtain a closed 2-cycle in $\mathcal O_{(X,a)}$ with
the desired properties. 

So it remains to check that we indeed have 
$\pi_0(\mathcal F_g)=\pi_1(\mathcal F_g)=\{0\}$ for each $g\in G$.
This is the content of the following lemma \ref{pi}.
\end{proof}

\begin{Lemma}\label{pi}
We have $\pi_0(\mathcal F_g)=\pi_1(\mathcal F_g)=\{0\}$ for all $g\in G$.
\end{Lemma}
\begin{proof}
Using the results of section \ref{coad} and denoting by $Stab_G(g)$
the stabilizer of $g$ with respect to $\tau$--twisted 
conjugation, we can write 
\begin{align*}
   \mathcal F_g = & \{(Y,a)\in\mathcal O_{(X,a)}~|~z_{(Y,a)}(2\pi)=g\}\\
   \cong& \{\gamma\in L(G,\tau)~|~\gamma(0)\in
       Stab_G(g)\}/Stab_{L(G,\tau)}(X,a)\,.
\end{align*}
Again using the results described in section \ref{coad}, 
we have $Stab_{L(G,\tau)}(X,a)\cong Stab_G(g)$ via the map
$\gamma\mapsto\gamma(0)$.
Let us write $\mathfrak G_g=\{\gamma\in L(G,\tau)~|~\gamma(0)\in
       Stab_G(g)\}$. Then we can use the long exact sequence 
$$
  \cdots\to\pi_1(Stab_G(g))
  \to\pi_1(\mathfrak G_g)\to\pi_1(\mathcal F_g)
   \to\pi_0(Stab_G(g))\to\pi_0(\mathfrak G_g)
   \to\pi_0(\mathcal F_g)\to 0
$$
to compute the homotopy groups. Indeed, it is easy to see that
$\mathfrak G_g$ is connected so that $\pi_0(\mathcal
F_g)=\{0\}$. Furthermore, since $G$ is simply connected and since
$Stab_G(g)$ is the fixed point set of an automorphism of $G$, we know
that $Stab_G(g)$ is connected. So if we can show that the injection
$\iota:Stab_G(g)\to \mathfrak G_g$ induces a surjection of fundamental
groups, we are done.

Let $\varphi: S^1\to\mathfrak G_g$, $\theta\mapsto\varphi_\theta$ 
be a loop in $\mathfrak G_g$. The map $\theta\mapsto
\varphi_\theta(0)$ defines a loop in $Stab_G(g)$. Obviously, this map
induces a surjection $\pi_1(\mathfrak G_g)\to\pi_1(Stab_G(g))$. 
It remains to show that the map is injective as well. So let
$\varphi$ and $\widetilde\varphi$ be two loops in $\mathfrak G_g$ which map
to the same element in $\pi_1(Stab_G(g))$. We have to show that
$\varphi$ and $\widetilde\varphi$ are homotopic. Let $f$ be the loop
in $Stab_G(g)$ defined via $f(\theta)=\varphi_\theta(0)$ and
accordingly $\widetilde f$.
Assume that $\Phi:[0,1]\times[0,2\pi]$
defines a homotopy from $f$ to $\widetilde f$. Then
$$
  \widetilde\Phi:(s,\theta)\mapsto
  \varphi_\theta\iota(\Phi(0,\theta))^{-1}\iota(\Phi(s,\theta))
$$ 
defines a
homotopy from $\varphi$ to a loop $\widehat\varphi$ in $\mathfrak G_g$
whose residual image $\widehat
f:\theta\mapsto\widehat\varphi_\theta(0)$ equals
$\widetilde f$. So from now on, we can assume that $f=\widetilde
f$. Finally, we can use the fact that the set $\{\gamma\in
L(G,\tau)~|~\gamma(0)=g_0\}$ is connected for each $g_0\in G$ to find
a homotopy from $\varphi$ to $\widetilde\varphi$. Indeed, we can view
the loops $\varphi$ and $\widetilde\varphi$ as sections in a fibration
over $S^1$, whose fiber at a point $\theta\in S^1$ is given by the set 
$\{\gamma\in L(G,\tau)~|~\gamma(0)=\varphi_\theta(0)\}$. One easily
checks that the
fibers are connected. Therefore, since $S^1$ is one-dimensional, the sections
$\varphi$ and $\widetilde\varphi$ are homotopic.
This finishes the proof.
\end{proof}


\end{document}